# Analysis of a Complex approximation to the Li-Keiper coefficients around the K Function

Danilo Merlini[1], Massimo Sala[2] and Nicoletta Sala[3]

[1] *CERFIM/ISSI, Locarno, Switzerland, e-mail: merlini@cerfim.ch*
[2] *Independent Researcher*
[3] *Institute for the Complexity Studies, Rome, Italy*

**Abstract**

We introduce a kind of "perturbation" for the Li-Keiper coefficients around the Koebe function (the K function) and establish a closed system of Equations for the Li-Keiper coefficients. We then check the correctness of some of the many possible solutions offered by the system, related to the discrete derivative of order n of a function.
We also report numerical finding which support our stability conjecture that the tiny part $\lambda_{\text{tiny}}(n)$ (the fluctuations around the trend) are bounded in absolute values by $\gamma \cdot n$, where $\gamma$ is the Euler-Mascheroni constant.

Paper's organization

1. Introduction.
2. The approximation and the K function.
3. The conjecture .
4. Numerical experiments.
5. Additional set of Equations and comparison for the first few Li-Keiper coefficients to 20 Digits.
6. Concluding remark.

**1. Introduction**

In the study of the Riemann Zeta function , in connection with the Riemann Hypothesis, the Li-Keiper coefficients play a fundamental role. It is in fact known that the non negativity property of the Li-
- Keiper coefficients is equivalent to the truth of the Riemann Hypothesis (RH).The,n-ten coefficient, usually indicated with $\lambda(n)$ is the coefficient of $z^n$ in the Taylor expansion around z=0 (s =1) of the log of the $\xi$ function written in the variable z = 1-1/s, i.e.

$$\log(\xi(1/(1-z))) = \sum_{n=1}^{\infty} z^n \cdot \lambda(n)/n$$

Usually one split the coefficients in two parts i.e. $\lambda(n) = \lambda_{trend}(n) + \lambda_{tiny}(n)$ i.e. the trend and the tiny part respectively defined as:

$\lambda_{trend}(n)$ = coefficient of $z^n$ in $\log(s.\Gamma(s/2)) = \log(1/(1-z).\Gamma((1/(1-z)/2))$
and $\lambda_{tiny}(n)$ = coefficient of $z^n$ in $\log((s-1) \cdot \zeta(s)) = \log((z/(1-z)) \cdot \zeta(1/(1-z)))$.
Calculations was made in terms of two sequences i.e. the Stieltjes constants $\gamma$'s and the other set of constants $\sigma$'s with coupled Equations among them. The trend behaves as $\lambda_{trend}(n) \sim (1/2).n.\log(n)+c.n$ while some numerical experiments report values of the tiny part $\lambda_{tiny}(n)$ up to some thousand of n , values which appear small with respect to the trend. There is also a gratuitous conjecture by us that $|\lambda_{tiny}(n)| < n^\varepsilon$ for all $\varepsilon > 0$. For some works on the subject, See [1, 2, 3, 4, 5].
Our aim in the following Sections is twofold: we first show how the tiny part of the Li-Keiper coefficients may be seen as a "perturbation" around an important function i.e. the Koebe function (the K function) which played an important role in the de Branges's proof of the Bieberbach conjecture; then we advance a conjecture which permits an extensivity behavior of the tiny part i.e. as $\gamma \cdot n$ where $\gamma$ is the Euler-
-Mascheroni constant; we then check some solutions of the general closed set of Equation we obtained ,-finally related - to the discrete derivative of a function.
We also report some numerical results which support our conjecture.

## 2. The approximation, the K function (the Koebe function), discrete derivative

We now introduce a kind of "perturbation" of the coefficients from those of the Koebe function , i.e. a perturbation around the Koebe function (See Appendix 1).We then consider (for values Re(s)~ 1) ,the function

$$f(s) = (1/\gamma) \cdot s \cdot (s-1), d/ds (\log[(s-1) \cdot \zeta(s)] ) \qquad (1)$$

where $\gamma = \lambda_{tiny}(1) = 0.577215....$ is the first coefficient of the "tiny" part of the expansion above , that is the Euler-Mascheroni constant [2] .
After the change of variable , $s = 1/(1-z)$ i.e. $z = 1-1/s$, here now $|z| < 1$ (the open unit disk), we have :

$$f(1/(1-z)) = (1/\gamma) \cdot z/(1-z)^2 \, (dz/ds) \cdot d/dz \left( \sum_{n=1}^{\infty} (\lambda_{tiny}(n)/n) \cdot z^n \right) =$$

$$= (1/\gamma) \cdot z \cdot \left( \sum_{n=1}^{\infty} n \cdot (\lambda_{tiny}(n)/n) \cdot z^{n-1} \right) = \sum_{n=1}^{\infty} (\lambda_{tiny}(n)/ \lambda_{tiny}(1)) \cdot z^n \qquad (2)$$

with $a_n = (\lambda_{tiny}(n)/ \lambda_{tiny}(1)$ , $n = 1,2,...$

the sequence of the Taylor coefficients $a_n$, ($a_0 = 0$, $a_1 = 1$); if the function $f(1/(1-z))$ above ,holomorphic in the interior of the unit disk i.e. $\pi$ for $|z| < 1$ where to be univalent then, the theorem of de Branges [7] would imply that $a_n < n$, i.e. our conjecture would be true since then $|\lambda_{tiny}(n)| < \gamma \cdot n$ and this will be a proof of the RH (Riemann Hypothesis) .
In fact the term $\gamma \cdot n$ would not be sufficient to eliminate the positivity property of the main term (the trend) for big values of n and related to the log of the Gamma function that for large n behaves as

$$\lambda_{trend} \sim (n/2) \cdot \log(n). \qquad (3)$$

(It is known that the trend behaves as $(n/2) \cdot \log(n) + c \cdot n$ for $n > n_0$ where c is the constant given by $(½) \cdot (\gamma - 1 - \log(2\pi)) = -1.13....$ [1, 3].
It should be said that the univalent property is difficult to be proven i.e. that for $z_1$ different from $z_2$ one has

f($z_1$) ≠ f($z_2$) (injection!). This appears as a necessary condition
in the de-Branges's Theorem [7]. See Appendix 2 for an illustration
where we limit ourselves to investigate the univalent condition in
two special cases.
Our conjecture [8] remains in any case independent of the de-Branges's
Theorem. We now proceed and consider the following strategy:
We reconsider the function f above :

$$f(1/(1-z)) = (1/\gamma).z/(1-z)^2 \cdot \Psi(z) = (1/\gamma).K(z).\Psi(z) \qquad (4)$$

where $\Psi(z) = (dz/ds) \cdot d/dz \left( \sum_{n=1}^{\infty} (\lambda_{tiny}(n)/n).z^n \right) =$

$= (1-z)^2 \cdot \Sigma \, n \cdot (\lambda_{tiny}(n)/n).z^{n-1} = (1-z)^2 \cdot \Sigma \, (\lambda_{tiny}(n)).z^{n-1} =$

$= \lambda_{tiny}(1) \cdot (1-2.z+z^2) + \lambda_{tiny}(2)) \cdot (1-2.z+z^2).z + \lambda_{tiny}(3)) \cdot (1-2.z+z^2).z^2 +$

$\ldots \lambda_{tiny}(n).(1-2.z+z^2).z^{n-1} = \lambda_{tiny}(1) + (\lambda_{tiny}(2) - 2.\lambda_{tiny}(1)).z^1 +$

$(\lambda_{tiny}(3) - 2.\lambda_{tiny}(2) + \lambda_{tiny}(1)).z^2 + \ldots\ldots (\lambda_{tiny}(n) - 2.\lambda_{tiny}(n-1) + \lambda_{tiny}(n-2)).z^n \ldots$

and where $K(z) = z/(1-z)^2$ above is the Koebe function (the K function)
i.e. $K(z) = z/(1-z)^2$ where in the Taylor expansion $a_0 = 0$, $a_1 = 1$ and
$a_n = n$, $n \geq 2$ [7], i.e.

$$K(z) = \sum_{n=1}^{\infty} n.z^n \,,\text{ which enters as "maximum" in the proof of}$$

de-Branges's Theorem. Notice here - as a comment - that if
in $\Psi(z)$ we consider for the $\lambda_{tiny}(n)$ the approximation given by

$$\lambda_{tiny}(n) - 2.\lambda_{tiny}(n-1) + \lambda_{tiny}(n-2) = 0 \qquad (5)$$

(coefficient of $z^n$), then $\Psi(z) = \gamma$ and $f(1/(1-z)) = K(z)$, i.e. it appears the
Koebe function; for a short discussion on the Koebe function, see
Appendix1; then with Eq. (5) above, we have

$$\sum_{n=1}^{\infty} (\lambda_{tiny}(n)/\lambda_{tiny}(1)).z^n \,, \text{ with } \lambda_{tiny}(n)/\lambda_{tiny}(1) = n \,. \qquad (6)$$

The function $(1/\gamma).\Psi(z)$ may be considered as a "perturbation" of the
Koebe function $K(z) = z + 2.z^2 + 3.z^3 + 4.z^4 +$ as the following
additional calculation shows:

$(1/\gamma).\Psi(z) = \quad 1 + [(\lambda_{tiny}(2)/\lambda_{tiny}(1) - 2].z +$

$\qquad + [(\lambda_{tiny}(3)/\lambda_{tiny}(1) - 2.\lambda_{tiny}(2)/\lambda_{tiny}(1) + 1].z^2 +$

$\qquad + [(\lambda_{tiny}(4)/\lambda_{tiny}(1) - 2.\lambda_{tiny}(3)/\lambda_{tiny}(1) + \lambda_{tiny}(2)/\lambda_{tiny}(1)].z^3 +$

(7)

and

$K(z).(1/\gamma).\Psi(z) =$
$= z + z^2.(2 + ((\lambda_{tiny}(2)/\lambda_{tiny}(1) - 2)) + z^3.[3 + 2.(\lambda_{tiny}(2)/\lambda_{tiny}(1) - 2) + ((\lambda_{tiny}(3)/\lambda_{tiny}(1) - 2.\lambda_{tiny}(2)/\lambda_{tiny}(1) + 1] +..$

$= z + z^2.(\lambda_{tiny}(2)/\lambda_{tiny}(1)) + z^3.(\lambda_{tiny}(3)/\lambda_{tiny}(1)) +...$  (8)

of course:

$(1/\gamma).z/(1-z)^2 \cdot \Psi(z) = z/(1-z)^2.(1/\lambda_{tiny}(1)).(1-z)^2.\Sigma (\lambda_{tiny}(n)).z^{n-1} =$

$\sum_{n=1}^{\infty} (\lambda_{tiny}(n)/\lambda_{tiny}(1)).z^n = \sum_{n=1}^{\infty} n.a_n.z^n \quad \text{where} \quad a_n = \lambda_{tiny}(n)/n.\lambda_{tiny}(1),$

that is using the Taylor expansion, we obtain:

$(1/\gamma).\log(z/(1-z).\zeta(1/(1-z)] = \Sigma \lambda_{tiny}(n)/n \lambda_{tiny}(1)..z^n =$

$= \gamma z + 0.483442.z^2 + 0.406898.z^3 + 0.343897.z^4 + 0.291653.z^5 +$
$\quad + 0.248049.z^6 +..$

$(1/\gamma).z/(1-z)^2 \cdot \Psi(z) =$

$= 0 + 1.z + 2.(0.483442/\gamma).z^2 + 3.(0.406898/\gamma).z^3 + 4.(0.343897/\gamma).z^4$
$+ 5.(0.291653/\gamma).z^5 + (0.248049/\gamma) z^6 + ... =$

$= 0 + 1.z + 2.(0.837).z^2 + 3.(0.704).z^3 + 4.(0.595).z^4 + 5.(0.505.)z^5 +$
$6.(0.429).z^6 + 7.(0.366).z^7 + 8.(0.312).z^8 + 9.(0.267).z^9 + 10.(0.229).z^{10}$

and the coefficients (0.837..) and so on are smaller than unity (the unity is reached only by the Koebe function as discussed above and we may now indicate our modified [8, 9] conjecture.

## 3. The conjecture

The series in z for the tiny fluctuations around the trend of the Li-Keiper coefficients is a perturbation of the series in z for the Koebe function and $a_n = \lambda_{tiny}(n)/n\,\lambda_{tiny}(1) < 1$, for $n \geq 1$.

By vanishing perturbation we have as derived above

$$\lambda_{tiny}(n) - \lambda_{tiny}(n-2) + \lambda_{tiny}(n-2) = 0 \qquad \lambda_{tiny}(0) = 0. \qquad (9)$$

Thus $\lambda_{tiny}(2) = 2.\,\lambda_{tiny}(1)$, $\lambda_{tiny}(3) = 3.\lambda_{tiny}(2) - 3.\lambda_{tiny}(1) = 3.\lambda_{tiny}(1)$
Thas is, $\lambda_{tiny}(n) = n.\lambda_{tiny}(1) = n.\gamma$ with $\gamma$ = the Euler-Mascheroni constant.

The above recurrence relation furnishes us an "approximation" of the $\lambda_{tiny}$'s at least for small n, i.e. by introducing now the true values obtained from the Taylor expansion (Approximation A).

$$\lambda_{tiny}(n) = 2.\,\lambda_{tiny}(n-1) - \lambda_{tiny}(n-2) \qquad (A)$$

we obtain, for the first values of n:

| n | $\lambda_{tiny}(n)/n.\lambda_{tiny}(1)$ | exact (3 decimals) |
|---|---|---|
| 2 | 1 | 0.837 |
| 3 | 0.783 | 0.704 |
| 4 | 0.638 | 0.595 |
| 5 | 0.530 | 0.505 |
| 6 | 0.444 | 0.429 |
| 7 | 0.375 | 0.366 |
| 8 | 0.318 | 0.312 |
| 9 | 0.271 | 0.267 |
| 10 | 0.231 | 0.229 |
| 11 | 0.199974 | 0.196782 |

Table1

The values of the first column appear as upper bounds of the exact one to 3 decimals given in the second column. For a numerical experiment in the (A)approximation (first order!) we refer to a previous work of us [9], where we have also studied the approximation (B) (second order discrete derivative) given by the relation similar to that of Eq.(A) above, i.e.

$$\lambda_{tiny}(n) - 3 \cdot \lambda_{tiny}(n-1) + 3 \cdot \lambda_{tiny}(n-2) + \lambda_{tiny}(n-3) = 0 \quad (B)$$

In this second approximation the corresponding Table similar to the above one is given below.

| n | $\lambda_{tiny}(n)/n.\lambda_{tiny}(n)$ | exact (3 decimals) |
|---|---|---|
| 2 | - | 0.837 |
| 3 | 0.675 | 0.704 |
| 4 | 0.579 | 0.595 |
| 5 | 0.496 | 0.505 |
| 6 | 0.424 | 0.429 |
| 7 | 0.362 | 0.366 |
| 8 | 0.310 | 0.312 |
| 9 | 0.266 | 0.267 |
| 10 | 0.228 | 0.229 |
| 11 | 0.196000 | 0.196782 |

Table 2

Here too, the values of the first column appear as lower bounds of the exact one to 3 decimals given in the second column.
At this moment it is tempting to consider further approximations as z-> 0 in the perturbation introduced above i.e.

$$\Psi(z) = (1-z)^m \cdot \sum_{n=1}^{\infty} n \cdot (\lambda_{tiny}(n)/n) \cdot z^{n-1}$$ i.e. with a power m instead

of the power 2 (Approx.(A)) or power 3 (Approx. (B)).

In fact approximations (A) and (B) suggest a more complex kind of "approximation", that is, the "tiny" fluctuations of order n, $\lambda_{tiny}(n)$ should contain more memory, thus depending from all previous tiny fluctuations $\{\lambda_{tiny}(k)\}$, k=1..n-1. (See below)

Then,
$$\lambda_{tiny}(n) = \varphi(\{\lambda_{tiny}(k)\}, k=1..n-1). \quad (10)$$

Our function $\varphi$ related to the discrete derivative turn out to be linear in all the $\lambda_{tiny}(k)$'s and is given by:

$$\lambda_{tiny}(n) = (-1)^{(n)} \cdot \sum_{k=1}^{n-1} (-1)^{(k-1)} \cdot \binom{n}{k} \cdot \lambda_{tiny}(k), \quad n \geq 1 \quad (D)$$

The first few terms of the sequence are:
- $\lambda_{tiny}(1) = \gamma = 0.577...$ (Euler-Mascheroni constant)
- $\lambda_{tiny}(2) = 2 \cdot \lambda_{tiny}(1) - 0$
- $\lambda_{tiny}(3) = 3 \cdot \lambda_{tiny}(2) - 3 \cdot \lambda_{tiny}(1) + 0$
- $\lambda_{tiny}(4) = 4 \cdot \lambda_{tiny}(3) - 6 \cdot \lambda_{tiny}(2) + 4 \cdot \lambda_{tiny}(1) - 0$
- $\lambda_{tiny}(5) = 5 \cdot \lambda_{tiny}(4) - 10 \cdot \lambda_{tiny}(3) + 10 \cdot \lambda_{tiny}(2) - 5 \cdot \lambda_{tiny}(1) + 0$
- $\lambda_{tiny}(6) = 6 \cdot \lambda_{tiny}(5) - 15 \cdot \lambda_{tiny}(4) + 20 \cdot \lambda_{tiny}(3) - 15 \cdot \lambda_{tiny}(2) + 6 \cdot \lambda_{tiny}(1) - 0$

.

$\lambda_{tiny}(n) = n \cdot \lambda_{tiny}(n-1) - .....$

## 4. Numerical experiments

In our numerical experiment with some digits, we have obtained the results, column (D) below for the values $\lambda_{tiny}(n)/n$ obtained with (D) after inserting the true values of the lambda's of lower index given by the exact sequence (C) (from "Tables"). After small values of n the agreement between column (D) and (C) is very satisfactory. Our values of (D) oscillate around the exact values of (C) (extremely tiny oscillations). As illustration for the reader we also give the plots of the two discrete functions below up to n =15, with the first 12 digits.

| n  | $\lambda_{tiny}(n)/n$ (D) | $\lambda_{tiny}(n)/n$ (C) |
|----|---------------------------|---------------------------|
| 2  | 0.5777215664901           | 0.483442548481            |
| 3  | 0.389669432061            | 0.406898976072            |
| 4  | 0.347584947674            | 0.343897032967            |
| 5  | 0.290748804461            | 0.291653700039            |
| 6  | 0.248290853736            | 0.248049721202            |
| 7  | 0.211388470880            | 0.211455834319            |
| 8  | 0.180626365551            | 0.180606968014            |
| 9  | 0.154505006842            | 0.154510711865            |
| 10 | 0.132382072505            | 0.132380368369            |
| 11 | 0.113585191797            | 0.113585706892            |
| 12 | 0.097616677732            | 0.097616520578            |
| 13 | 0.084055437619            | 0.084055485931            |
| 14 | 0.072557830130            | 0.078557815185            |
| 15 | 0.062835890090            | 0.062835894739            |

Table 3

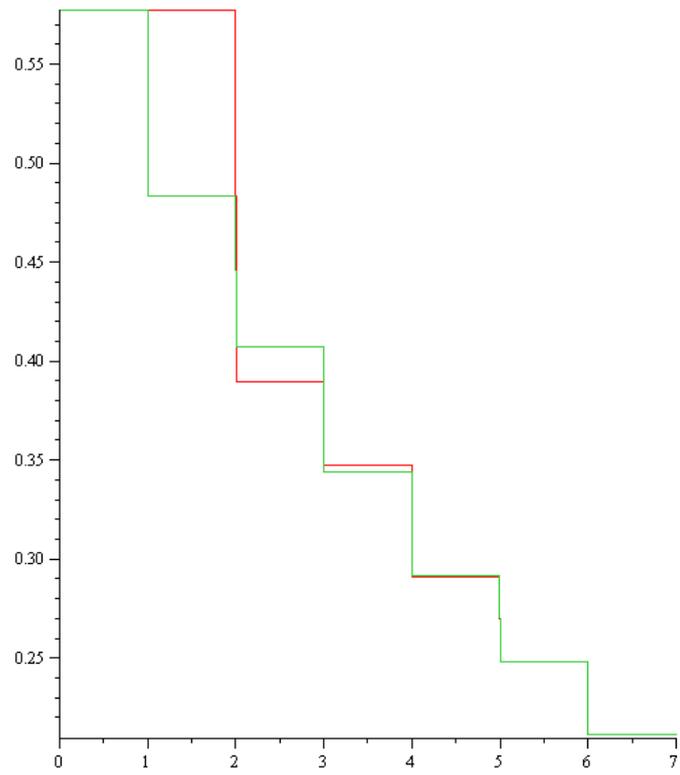

Fig. 1 (D) in red and ( C) in green, in the range n = [0..7]. The values of (D) are alternating to the values of ( C ).

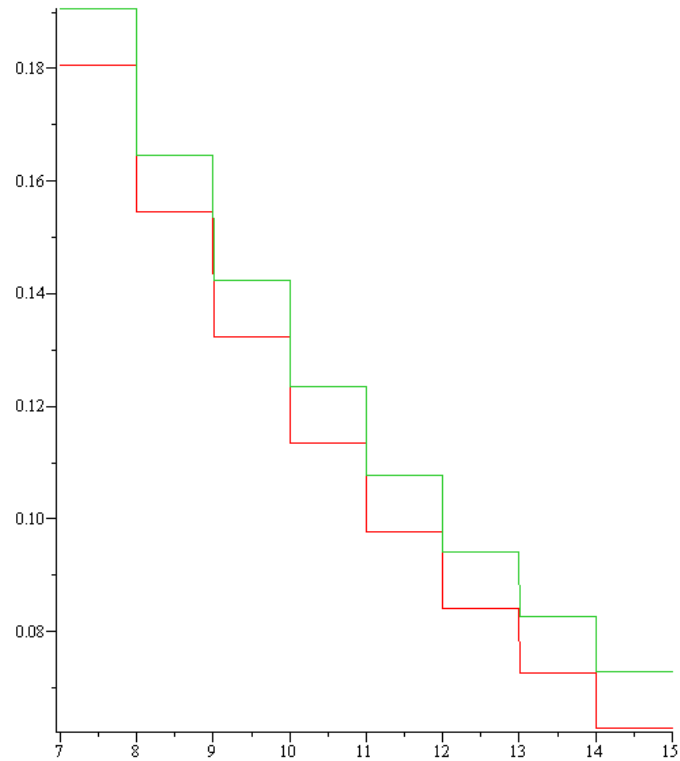

Fig. 2 The Plots of (D) in red and (C) in green (increased by 0.01),
in order that their "difference" be visible, in the range
n= [7..15]. (at n=15 the error is only of the order $4.649 \cdot 10^{-9}$, with
a relative error of ~ $7.4 \cdot 10^{-8}$).

We now look at the above Formula for the trend, i.e. in order to find analogous "approximations" for $\lambda_{trend}(n)/n$ (E) to be compared with the true $\lambda_{trend}(n)/n$ (F). For (E) we apply the same Formula as (D) above where instead of $\lambda_{tiny}(n)/n$ we now use $\lambda_{trend}(n)/n$ from known Tables. The computations was carried out with 30 digit; we give here the results with 12 digits.

For the trend too, we notice the alternating property between (E) and (F). Here too, we verify that the error at n=15 is of the order $4.65 \cdot 10^{-9}$, with a relative error of $1.699 \cdot 10^{-8}$.

| n | $\lambda_{trend}(n)/n$ (E) | $\lambda_{trend}(n)/n$ (F) |
|---|---|---|
| 1 | - 0.554119955935 | - 0.554119955935 |
| 2 | - 0.554119955935 | -0.437269680867 |
| 3 | - 0.320419405799 | - 0.337686002554 |
| 4 | - 0.255368920995 | - 0.251699413094 |
| 5 | - 0.175640404588 | - 0.176545157147 |
| 6 | - 0.110363231928 | - 0.110122052488 |
| 7 | - 0.050751310573 | - 0.050818674666 |
| 8 | + 0.002593094250 | + 0.002612991628 |
| 9 | 0.051152331863 | 0.051146626843 |
| 10 | 0.095551863813 | 0.095553567949 |
| 11 | 0.136447611677 | 0.136447096581 |
| 12 | 0.174321265653 | 0.174321422807 |
| 13 | 0.209578412983 | 0.209578364671 |
| 14 | 0.242547718346 | 0.242547733292 |
| 15 | 0.273502734711 | 0.273502730062 |

Table 4

Below we present the plots of the discrete functions (E) and (F).
and finally the plot of the complete $\lambda(n)/n$ , the discrete function
(G) = (D) + (E)  with the true function (H) = (C) + (F).
Notice that if we multiply the values of (H) by n, i.e. for n=1 by 1,
for n=2 by 2 and so on we obtain the values of the last column in the
table of Ref.[1] for the $\lambda$'s. Here the values of the Tables are X/n not X
(X= $\lambda_{trend}(n)$ , $\lambda_{trend}(n)$ or $\lambda(n)$)). For Example, for n= 9, from our Table we
have: (0.1545107 +0.0511446626863).9 = 1.85.91... as in the Ref [1].

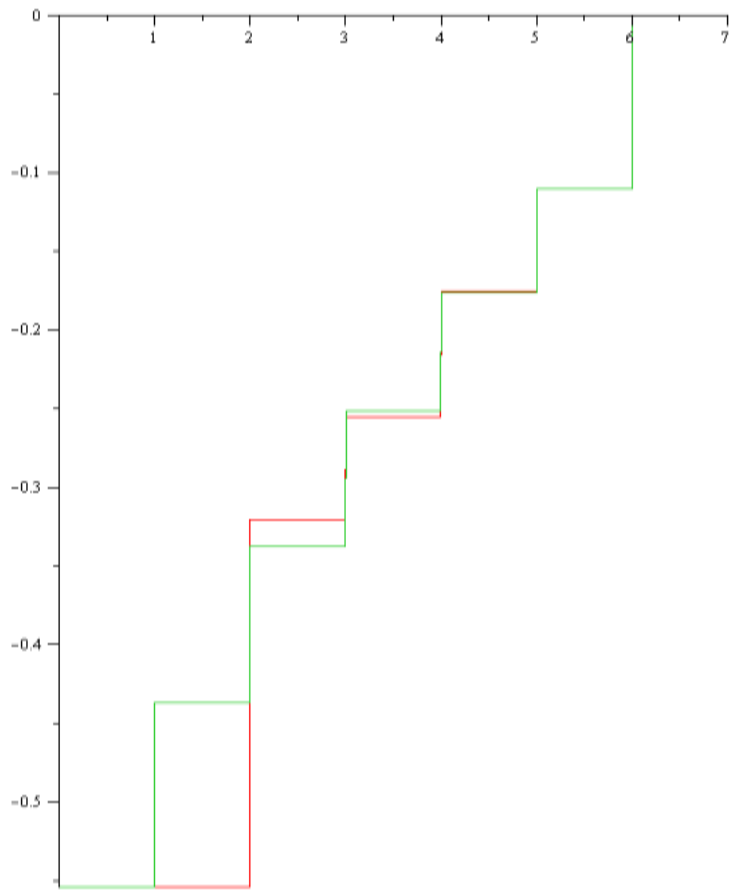

Fig 3 (E) in red, (F) (the true function!) in green , in the range n= [0, 7]  for the function  $\lambda_{trend}(n)/n$ .

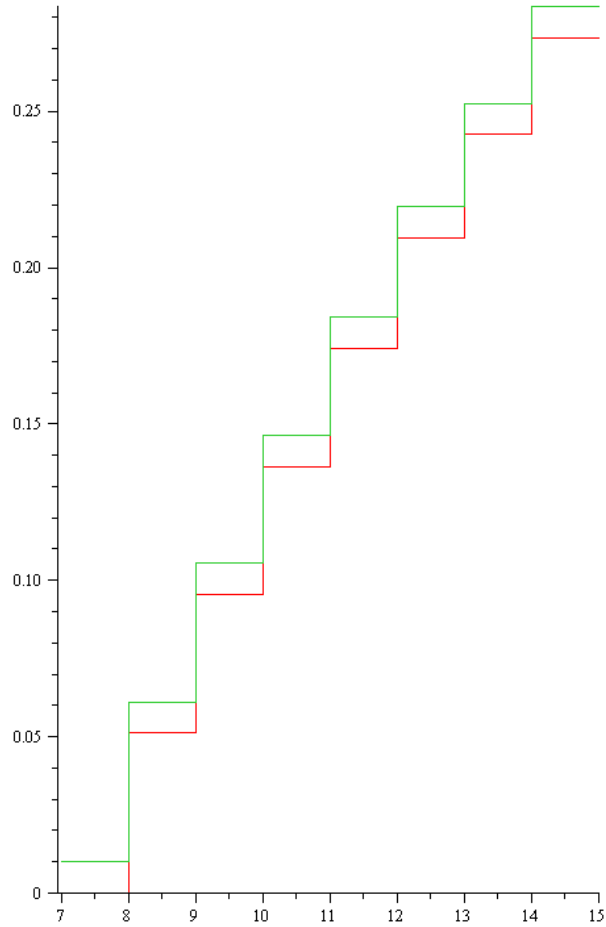

Fig. 4 The Plots of (E) in red and of (F) -increased by 0.01- ,in order that their "difference" be visible, in the range n = = [7..15].

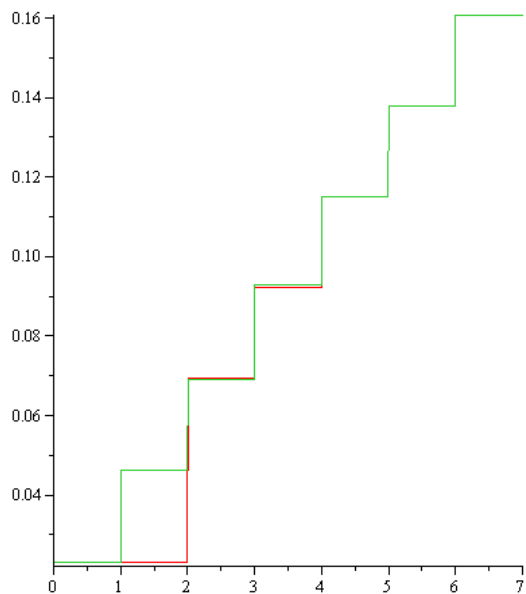

Fig. 5 Plot of the discrete functions (G) in red and (H) in green (the true function) in the range n = [ 0,7].

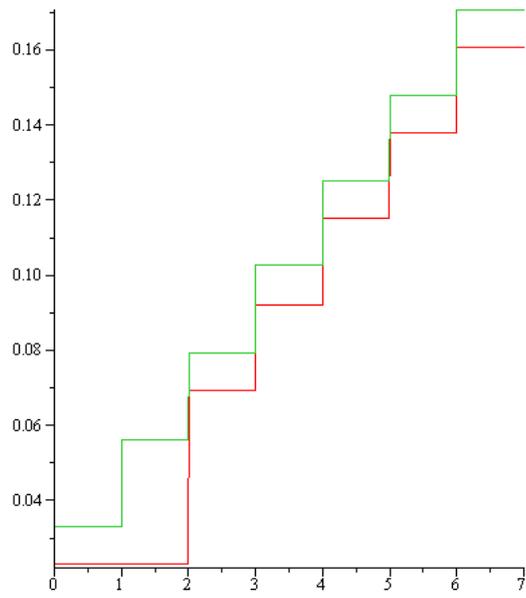

Fig. 6 Plot of the discrete functions (G) in red and and (H) (increased by 0.01), in green, in the range n= [0,7 ].

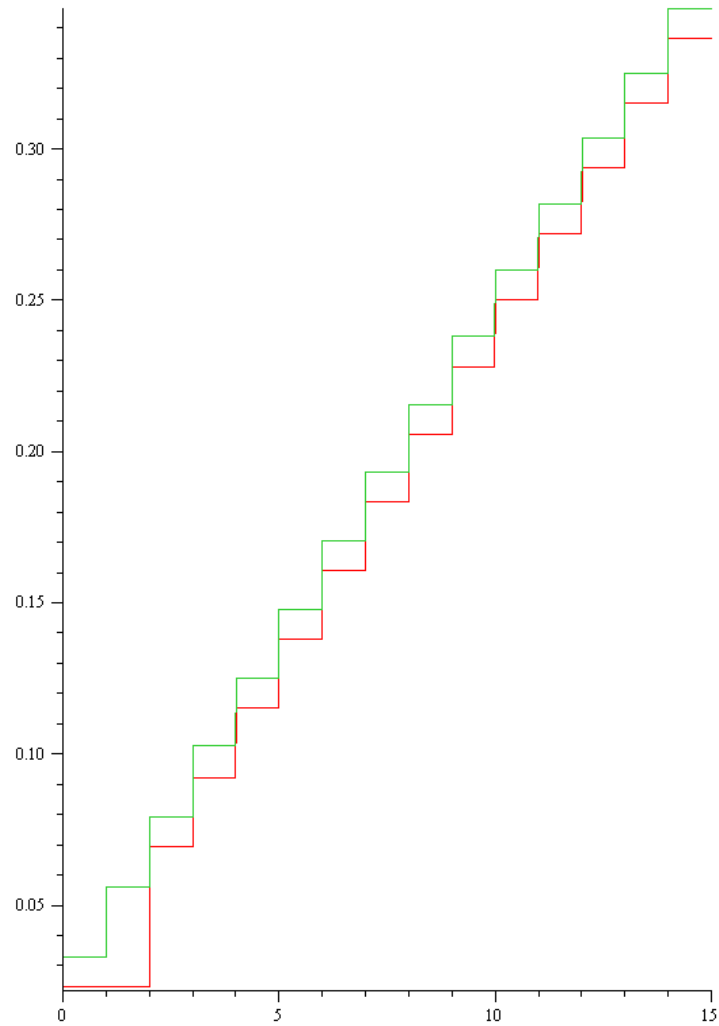

Fig. 7 Plot of (G) in red and (H) in green in the range n= [0..15] .
(We still increased (H) by 0.01 in order that the two plots be "visible").

Before our conclusion we may add the following: For $\lambda_{tiny}(n)/n$ and $\lambda_{trend}(n)/n$ , we have seen that both converge fast to the true values - using our approximation- which involves the whole history of the $\lambda$'s. (discrete derivatives).

For the complete $\lambda(n)/n$ or $\lambda(n)$ which is the sum of the two contributions $\lambda_{trend}(n)/n + \lambda_{tiny}(n)/n$ resp. $\lambda_{trend}(n) + \lambda_{ttiny}(n)$ we may easily relate the kind of approximations we have considered to the general discrete derivative along the above lines. In fact from the definition of the Li-Keiper coefficients we have:

$$\lambda(n) = \sum_{\rho} (1-(1- 1/\rho)^n) \qquad (11)$$

where the sum is over all nontrivial zeros of the Zeta function in the critical strip, i.e. for $0 < \sigma < 1$ ($s = \sigma + i.t$). Now, we will derive a system of closed Equations for the Li-Keiper coefficients, using (11). The system may be applied to the trend and to the tiny part since such a system is obtained directly from the definition of the coefficients but it is then recognized that the central point is the discrete derivative of functions. We have:

$$\lambda(1) = \sum_\rho (1-(1- 1/\rho )) = \sum_\rho (1/\rho) = (1+\gamma-(1/2).\log(4.\pi))= 0.0230957...$$

$$\lambda(2) = \sum_\rho (1-(1- 1/\rho))^2 = \sum_\rho (2/\rho - 1/\rho^2)$$

$$\lambda(2) = 2.(\sum_\rho(1-(1- 1/\rho)) - \sum_\rho (1/\rho^2) = 2.\lambda(1) - \sum_\rho (1/\rho^2)$$

$$\lambda(3) = \sum_\rho(1- (1- 1/\rho)^3) = \sum_\rho (3/\rho - 3/\rho^2 + 1/\rho^3) =$$

$$= - \sum_\rho (3/\rho) + \sum_\rho ( 6/\rho - 3/\rho^2 + 1/\rho^3) =$$

$$= 3.( \sum_\rho (1- (1-1/\rho)^2) - 3.\sum_\rho (1- (1-1/\rho)) + \sum_\rho 1/\rho^3 =$$

$$= 3.\lambda(2) - 3.\lambda(1) + \sum_\rho 1/\rho^3$$

In the same way we obtain:

$$\lambda(4) = 4.\lambda(3) - 6.\lambda(2) + 4.\lambda(1) - \sum_\rho 1/\rho^4$$

with the appearance of the binomial coefficients; we thus have in general, that for $n > 1$:

$$\lambda(n) = \sum_{k=1}^{n-1}(-1)^{(k-n+1)} \cdot \binom{n}{k} \cdot \lambda(k) + \sum_\rho (-1)^{(n-1)} \cdot \left(\frac{1}{\rho^n}\right)$$

$$\lambda(n) = \sum_{k=1}^{n-1}(-1)^{(k-n+1)} \cdot \binom{n}{k} \cdot \lambda(k) + \Delta \qquad (12)$$

with $\Delta = \sum_\rho (-1)^{(n-1)} .(1/\rho^n)$.

In fact, from the definition we have:

$$* \sum_{k=1}^{n-1}(-1)^{(k-n+1)} \cdot \binom{n}{k} \cdot \sum_\rho \left(1 - \left(1 - \frac{1}{\rho}\right)^k\right) + \sum_\rho (-1)^{(n-1)} \left(\frac{1}{\rho^n}\right) =$$

$$= \sum_{k=0}^{n}(-1)^{(k-n+1)} \cdot \binom{n}{k} \cdot \sum_\rho \left(1 - \left(1 - \frac{1}{\rho}\right)^k\right) - 0 + \sum_\rho \left(1 - \left(1 - \frac{1}{\rho}\right)^n\right) +$$

$$+ \sum_\rho (-1)^{(n-1)} \cdot \left(\frac{1}{\rho^n}\right) = \sum_{k=0}^{n}(-1)^{(k-n+1)} \cdot \binom{n}{k} \cdot \sum_\rho \left(-\left(1 - \frac{1}{\rho}\right)^k\right) +$$

$$+ \sum_\rho (-1)^{(n-1)} \cdot \left(\frac{1}{\rho^n}\right) + \lambda(n) = \lambda(n) + \sum_{k=0}^{n}(-1)^{(n-1)} \cdot \left(\frac{1}{\rho^n}\right) +$$

$$(-1)^{(n-1)} \cdot \sum_{k=0}^{n}(-1)^{(k-n+1)} \cdot \binom{n}{k} \cdot \left(1 - \frac{1}{\rho}\right)^k =$$

$$= \lambda(n) + \sum_\rho (-1)^{(n-1)} \cdot \left(\frac{1}{\rho^n}\right) - (-1)^{(n-1)} \cdot \left(1 - \left(1 - \frac{1}{\rho}\right)\right)^n = \lambda(n)$$

We now look at the kind of approximation we have if we neglect
Δ. As an example, for the sequence $\lambda(n) = 2.\lambda(n-1) - \lambda(n-2)$ the
characteristic Equation of the sequence is $x^2 = 2.x - 1$, with the solution
x=1. Thus $\lambda(n) = \alpha.n + \beta$, i.e. an arithmetic progression.
This holds in the general case as it should be and as illustration
we check that in the above linear approximation we obtain the
solution $\lambda(n) = \alpha.n + \beta$ (among others of course) We have:

$$\lambda(n) = \sum_{k=1}^{n-1}(-1)^{(k-n+1)} \cdot \binom{n}{k} \cdot \lambda(k) =$$

$$= \sum_{k=1}^{n-1}(-1)^{(k-n+1)} \cdot \binom{n}{k} \cdot (\alpha \cdot k + \beta) \cdot \beta \cdot (-1)^{(n-1)} \cdot \sum_{k=1}^{n-1}(-1)^{(k)} \cdot \binom{n}{k} =$$

$$= \beta \cdot (-1)^{(n-1)} \cdot \sum_{k=1}^{n-1}(-1)^{(k)} \cdot \binom{n}{k} - \beta \cdot (-1)^{(n-1)} - \beta \cdot (-1)^{(n-1)} \cdot (-1)^n$$

$$= 0 - \beta \cdot ((-1)^{(n-1)} + (-1)^{(2.n-1)}), \text{ equal to } 0. \beta \text{ if n is odd and to } 2. \beta$$

if n is even. Thus $\beta = 0$. Then since

$$\sum_{k=1}^{n-1}(-1)^{(k)} \cdot \binom{n}{k} \cdot k = n \cdot (-1)^{(n-1)}$$

we obtain $\alpha.n.(-1)^{(2n-2)} = \alpha.n$, $\rightarrow \lambda(n) = \alpha.n$.

We now compute an upper bound on

$$\Delta(n) = \sum_{\rho}(-1)^{(n-1)}.(1/\rho^n). \qquad (13)$$

Since the non trivial zeros are of the form $\rho = \sigma + i.t$ and that they may be paired as usually, then:

$$|\Delta(n)| < \sum_{\rho}[1/(\sigma^2 + t^2)]^n \sim \int_{t_0}^{\infty}(1/(2.\pi)).\log(t/2.\pi).(1/t^2)^n \cdot dt.$$

Keeping in mind that the density of the nontrivial zeros is essentially given by $dN = (1/(2.\pi)).\log(t/(2.\pi)).dt$, that an upper bound is given by setting $\sigma = 0$ and remembering that the value of $t_0$ of the first zero is $t = 14.134725..$, we have:

$$|\Delta(n)| < (1/(2.\pi)).(1/14)^{(2.n-1)}.(1/(2.n-1)).(\log(14/2.\pi)+(1/(2n-1)))$$

essentially, $|\Delta(n)| < c.(1/(14)^{2.n})$ with c a constant depending on n, thus small for big n. For n = 5, we have an amount of $(1/14^{10}) \sim 10^{-12}$. Notice that the above computation do not assume that the Riemann Hypothesis is true. Of interest is to check the correctness of some other solutions of our system of Equations with the binomial structure. We now check the correctness of the important function c.n.log(n) with c a free constant.

Let $\lambda(n) = n.\log(n)$, i.e assume $\lambda(k) = k.\log(k)$; then we have to compute $\varphi_1$, i.e.

$$\varphi_1 = \sum_{k=1}^{n-1}(-1)^k \cdot \binom{n}{k} \cdot (k \cdot \log(k)) \qquad (14)$$

to be compared with $\varphi_2 = (-1)^{(n-1)} \cdot n \cdot \log(n)$. Below we give the Table of the two functions in the range [0.. 31] of n, calculated up to

20 Digits. The difference becomes smaller as n increases, and $(\varphi_1 - \varphi_2)/n$ is vanishing as $n \to \infty$. Thus our function is:

$\varphi_1(n) = \varphi_2(n) \pm \delta(n) = (-1)^{(n-1)} \cdot n \cdot \log(n) \pm \delta(n)$, $|\delta(n)| < 0.30$, for $n > 31$.

which shows the correctness of the function $\varphi_1(n)$ if compared with $\varphi_2$.
We have thus controlled two families of solutions of the linear set of Equations i.e. of the form $\alpha \cdot n + \beta \cdot n \cdot \log(n)$.
For $\beta = \frac{1}{2}$ we recover the correct behavior for (the trend) and as $n \to \infty$, $\lambda(n) \sim (\frac{1}{2}) \cdot n \cdot \log(n)$ as the dominant term.
A point of interest is then still to check other interesting solutions of the form $\varphi(n) = n^\nu$, $\nu > 0$, also solutions of the form $\varphi \sim \log(n)$.
Below we give the Table of $\varphi_1$ and of $\varphi_2$ above.

| n | $\varphi_1$ | $\varphi_2$ |
|---|---|---|
| 1 | 0. | 0. |
| 2 | 0. | -1.386294361119 |
| 3 | 4.158883083359 | 3.295836866004 |
| 4 | -4.865581297297 | -5.545177444479 |
| 5 | 8.630462173553 | 8.047189562170 |
| 6 | -10.227797609117 | -10.750556815368 |
| 7 | 14.102018732148 | 13.621371043387 |
| 8 | -16.186194269648 | -16.635553233343 |
| 9 | 20.199979356843 | 19.773021196025 |
| 10 | -22.620534931015 | -23.025850929940 |
| 11 | 26.765919774058 | 26.376848000782 |
| 12 | -29.443522577892 | -29.818879797456 |
| 13 | 33.707924745111 | 33.344341646999 |
| 14 | -36.593468290124 | -36.946802614613 |
| 15 | 40.965062362406 | 40.620753016533 |
| 16 | -44.025136211015 | -44.361419558364 |
| 17 | 48.493711825837 | 48.164626848955 |
| 18 | -51.704110495157 | -52.026691642130 |
| 19 | 56.261007350245 | 55.944340604162 |
| 20 | -59.603387920029 | -59.914645471079 |
| 21 | 64.241256392224 | 63.934971192191 |
| 22 | -67.701240409001 | -68.002933973882 |
| 23 | 72.413803051248 | 72.116366966370 |
| 24 | -75.979818152899 | -76.273291928350 |
| 25 | 80.761669341091 | 80.471895621705 |
| 26 | -84.424202080048 | -84.710509988558 |
| 27 | 89.290064772599 | 88.987595382116 |
| 28 | -93.021739955872 | -93.301726284905 |
| 29 | 97.928670979204 | 97.651579069607 |
| 30 | -101.761680256697 | -102.035921449864 |
| 31 | 106.725360482290 | 106.453603339039 |
| 32 | -110.634257921228 | -110.903548889591 |

Table 5

We also add the Plots of the two functions in the interval of n = [0..20].

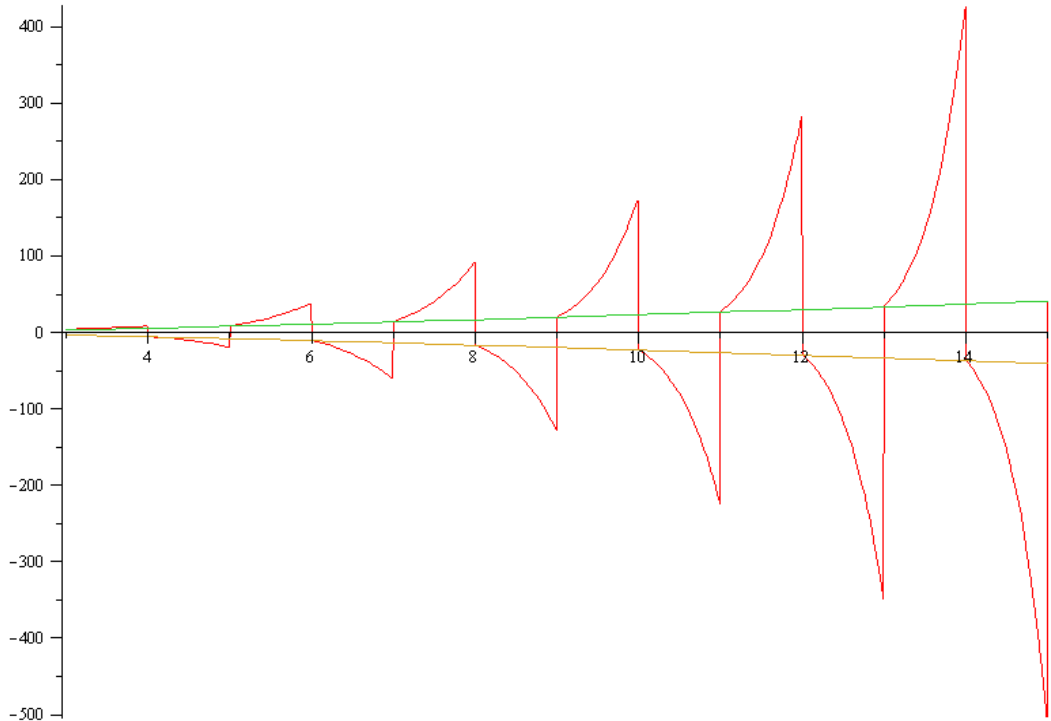

Fig. 8 The plot of $\varphi_1$ (in red) ,of $|\varphi_2|$ (in green) and $-|\varphi_2|$ in yellow.

As a last numerical experiment we give explicitly the formula where the fluctuations around the trend are given by our more plausible conjecture i.e. that the fluctuations are at most linear with n around the trend i.e. O(n), in connection with some known numerical computations reported in [1, 3 ,4 ]. The Formula is:

$$\lambda(n) = \sum_{k=1}^{n-1} (-1)^k \cdot \binom{n}{k} \cdot ((1/2) \cdot k \cdot \log(k) + c \cdot k + -\gamma \cdot k) \quad (15)$$

where $c = (½) \cdot (1+\gamma -\log(2\pi)) = -1.13...$ [ 1 ] and $\gamma = 0.577..$ is the Euler-Mascheroni constant [2].

## 5. An additional set of closed Equations and comparison for the first few Li-Keiper coefficients to 20 Digits.

The additional set of Equations for the lambda's of interest here( have been derived in Ref [10]. The set is different from our (Eq.(10)) and involves the binomial coefficients (2n,n-k) instead of (n,k) .as in Eq.(10), that is[10]:

$$\lambda_n = -n \cdot \sum_{j=1}^{n} (-1)^j (1/j) \cdot \binom{n+j-1}{2j-1} \cdot Z(j) \qquad (16)$$

where $Z(\sigma) = \sum_{\rho} (1/x_k^{\sigma})$ Re($\sigma$) > 1/2 and $x_k = \rho \cdot (1-\rho) = \frac{1}{4} + \tau_k^2$
Re($\tau_k$) > 0 , $\rho = \frac{1}{2} +- i \cdot \tau_k$ , k=1... $\infty$ , with $\rho$ any nontrivial zero of Zeta in the critical strip.

The inversion of the above Equation reads:

$$\sum_{k=0}^{n} (-1)^k \cdot \binom{2n}{n-k} \cdot \lambda_k = Z(n) \quad . \qquad (17)$$

Following these formula we may now show that Z(j) , is small as a function of j. In fact it may be shown that $|Z(j)| < (1/14.134..)^{(2 \cdot j-1)}$ where $t_1 = 14.134...$ is the height of the first zero of Zeta. We now introduce a new approximation following the above relations which clearly amount to set Z(j) = 0 ,in order to have a closet set of Equations for the Li-Keiper coefficients, called Approximation 2 below and given using Eq. (17) ,by:

$$\lambda_n = (-1)^{(n)} \cdot \sum_{k=1}^{n-1} (-1)^k \cdot \binom{2n}{n-k} \cdot \lambda_k \qquad (18)$$

The first few terms are:

$\lambda_2 = 4 \cdot \lambda_1$ .
$\lambda_3 = 6 \cdot \lambda_2 - 15 \cdot \lambda_1$
$\lambda_4 = 8 \cdot \lambda_3 - 28 \cdot \lambda_2 + 56 \cdot \lambda_1$ and so on instead of our previous

Approximation 1 (See Eq. (10)) given by:
$\lambda_2 = 2 \cdot \lambda_1$ .
$\lambda_3 = 3 \cdot \lambda_2 - 3 \cdot \lambda_1$
$\lambda_4 = 4 \cdot \lambda_3 - 6 \cdot \lambda_2 + 4 \cdot \lambda_1$ and so on, Approximation 1.

The numerical experiment we carried out is up to n=7, but from known numerical advanced calculations, we take values up to 20 Digits. [11].

Below we give the table of the first 7 Li-Keiper coefficients exact to 20 digits obtained from Approximations 1 (Eq.(10) and Approximation 2, (Eq.(18) together with the true values taken from advanced computations (Tables in Ref [ [ 11 ] ) up to n=7 with 20 digits.
(Approximation 1=$A_1$, Approximation 2 =$A_2$ ).

$\lambda_1 =$     0.02309570896612103381
$\lambda_2 =$     0.09234573522804667038
$\lambda_3 =$     0.20763389205543248037
$\lambda_4 =$     0.36879047949224163856
$\lambda_5 =$     0.57554271446117745240
$\lambda_6 =$     0.82756601228237929740
$\lambda_7 =$     1.12446011757095949058

$\lambda_2 (A_2) =$ **0.0923**828358644841352
$\lambda_2 (A_1) =$ 0.0461914179322420676

$\lambda_3 (A_2) =$ **0.20763**877687646451528
$\lambda_3 (A_1) =$ **0.2077**5007878577690974

$\lambda_4 (A_2) =$ **0.368790**48015206955248
$\lambda_4 (A_1) =$ **0.3688**6410671350332808

$\lambda_5 (A_2) =$ **0.5755427144**5798356265
$\lambda_5 (A_1) =$ **0.57554**199936782168970

$\lambda_6 (A_2) =$ **0.8275660122823**95050
$\lambda_6 (A_1) =$ **0.82756**5730084596235

$\lambda_7 (A_2) =$ **1.12446011757095943**661
$\lambda_7 (A_1) =$ **1.1244601**2214515063973 .

From the above numerical results we notice the excellent values with Approximation 2 ($A_2$) and good values too with Approximation 1($A_1$).

Remark:

Approximation 1 follows from a general behavior of the a discrete derivative of a function , while Approximation 2 follows from a specific treatment using the structure of the $\xi$ function , See Ref [10] . Approximation 2 has been derived using the structure of the zeros of the $\xi$ function while Approximation1- as we discussed - it is more generic but also valid for the trend and other functions.

The two Approximation have in common the following situation:
If instead to insert the true values of the $\lambda$'s in the Equations we use only the initial condition, we have for Approximation 1:

$\lambda_2(A_1) = 2.\lambda_1$
$\lambda_3(A_1) = 3.\lambda_2 - 3.\lambda_1 = 3.\lambda_1$ .
$\lambda_4(A_1) = 4.\lambda_3 - 6.\lambda_2 +4.\lambda_1 = 4.\lambda_1$ and so on, a linear behavior i.e.

$$\lambda_n(A_1) = n \cdot \lambda_1$$

for Approximation 2 we have instead:

$\lambda_2(A_2) = 4.\lambda_1$
$\lambda_3(A_2) = 6.\lambda_2 - 15.\lambda_1 = 9.\lambda_1$
$\lambda_4(A_2) = 8.\lambda_3 - 28\lambda_2 +56.\lambda_1 = (72 - 112 + 56).\lambda_1 = 16.\lambda_1$ and so on, a

quadratic behavior
$$\lambda_n(A_2) = n^2 \cdot \lambda_1 \ .$$

On the other hand if for $(A_1)$ we assume the initial condition $\lambda_2(A_1) = 4.\lambda_1$

we obtain:

$\lambda_3(A_1) = 3.\lambda_2 - 3.\lambda_1 = (12-3).\lambda_1 = 9 \cdot \lambda_1$ .
$\lambda_4(A_1) = (4 \cdot 9 - 6 \cdot 4 + 4) \cdot \lambda_1 = 16 \cdot \lambda_1$
$\lambda_5(A_1) = 5.\lambda_4 - 10 \cdot \lambda_3 + 10 \cdot \lambda_2 - 5 \cdot \lambda_1 = (80 - 90 + 40 - 5) \cdot \lambda_1 = 25 \cdot \lambda_1 = 5^2 \cdot \lambda_1$

and so on , we obtain a quadratic law too. This illustrate the influence of the initial condition for such a "process " .

Now, if we take for $A_1$, the initial condition :

$\lambda_2(A_1) = 3 \cdot \lambda_1$ instead of $\lambda_2(A_1) = 2 \cdot \lambda_1$ we obtain the sequence:
$\lambda_3(A_1) = (3 \cdot 3 - 3) \cdot \lambda_1 = (3+3)\lambda_1$
$\lambda_4(A_1) = 4 \cdot \lambda_3 - 6 \cdot \lambda_2 + 4 \cdot \lambda_1 = (24 - 18 + 4) \cdot \lambda_1 = 10 \cdot \lambda_1 = (3+3+4) \cdot \lambda_1$
$\lambda_5(A_1) = 5 \cdot \lambda_4 - 10 \cdot \lambda_3 + 10 \cdot \lambda_2 - 5 \cdot \lambda_1 = (5 \cdot 10 - 6 \cdot 10 + 10 \cdot 3 - 5) \cdot \lambda_1 = 15 \cdot \lambda_1$
$\quad = (3+3+4+5) \cdot \lambda_1$
$\lambda_6(A) = 6 \cdot \lambda_5 - 15 \cdot \lambda_4 + 20 \cdot \lambda_3 - 15 \cdot \lambda_2 + 6 \cdot \lambda_1 =$
$\quad = (6 \cdot 15 - 15 \cdot 10 + 20 \cdot 6 - 15 \cdot 3 + 6) \cdot \lambda_1 = 21 \cdot \lambda_1 =$
$\quad = (3+3+4+5+6) \cdot \lambda_1$, and so on and :

$$\lambda_n(A) = (3 + (3+4+5+6+\ldots n)) = n \cdot (n+1)/2$$

also a quadratic law.

Finally, for the general initial condition i.e. $\lambda_2(A_1) = c \cdot \lambda_1(A_1) = c \cdot \lambda_1 = c \cdot 0.0230957\ldots$ , for the approximation 1 we obtain:

$$\lambda_n(A_1) = [(c/2 - 1) \cdot n \cdot (n-1)] + n ] \cdot \lambda_1$$

i.e. and a linear law for $c=2$ and always a quadratic law for $c>2$.
Of course, inserting in $(A_1)$ and in $(A_2)$ the true values, the terms in $n \cdot \log(n)$ and in $n$ emerge from both set of Equations.

## 6. Concluding remark

In the first part of our work we have presented some analysis for the series in $z = 1 - 1/s$ of the tiny part of the Li-Keiper coefficients which enabled us to formulate the series for the tiny part as a "perturbation" around the Koebe function $z/(1-z)^2$, a function which appears as a "maximum" in the de Branges's Theorem for the Bieberbach conjecture.
In the second part, the aim was to find a good approximation for the complete Li-Keiper coefficients directly from the definition of them. The finding is a set of closed Equations where the coefficient $\lambda(n)$ contain the whole "history" of those of lower order. The closed system holds also for the trend as well for the tiny part since it may finally be seen that it is described by a formula concerning the discrete derivative of a function [6].
In connection with advanced numerical computations and in analogy with some kind of stability bounds appearing in statistical mechanics our present conjecture is finally reformulated as follow:

Conjecture : the tiny part of the Li-Keiper coefficients that is $\lambda_{tiny}(n)$ verify the inequality

$$| \lambda_{tiny}(n)/(n \cdot \lambda_{tiny}(1) ) | = | \lambda_{tiny}(n)/(n \cdot \gamma)| \leq 1$$

for any n. The above absolute value is "decreasing" with n and no value of n contradicting the above inequality has been found in a numerical experiment as reported below ( we have read approximately some of the numerical values of $\lambda_{tiny}(n)$ from [1, 3, 4].

| n | $\lambda_{tiny}(n)/(n \cdot \gamma)$ |
|---|---|
| 1 | 1 |
| 5 | 0.5052 |
| 8 | 0.3128 |
| 10 | 0.2293 |
| 100 | 0.0108 |
| 300 | 0.0277 |
| 800 | -0.0138 |
| 900 | 0.0173 |
| 1000 | 0.0030 |
| 1200 | 0.0083 |
| 1550 | - 0.0044 |
| 1800 | 0.0053 |
| 2000 | 0.0096 |
| 2440 | -0.0045 |
| 2650 | - 0.0036 |
| 3000 | - 0.0040 |
| 3240 | - 0.0069 |
| 3580 | 0.00387 |
| 4000 | - 0.0030. |
| 4860 | -0.0034 |
| 5100 | 0.0065 |
| 5500 | 0.0046 |
| 5700 | - 0.0035 |
| 5900 | - 0.00367 |
| 6200 | 0.00237 |
| 7000 | 0.00254 |
| 8000 | 0.00346 |
| 9630 | -0.000539 |
| 10000 | -0.000727 |

# References

1. K. Maslanka; "Effective method of computing Li's coefficients and their properties", ArX3iv:math/0402168v5 (math NT).

2. Yu.V. Matiyasevich: "Yet Another Representation for Reciprocals of the Nontrivial Zeros of the Riemann Zeta Function",Mat.Zametki,97:3(2015), 476-479.

3. J.B. Keiper: "Power Series Expansions of Riemann's function": Mathematics ,Volume 58 ,number 198,April 1992,765-773

4. F. Johansson: "Rigorous high-precision computation of the Hurwitz zeta function and its derivatives", ArXiv:1309.2877v1 (cs.SC) (2013).

5. D. Merlini, M. Sala, and N. Sala: "Fluctuation around the Gamma function and a Conjecture", IOSR International Journal of Mathematics, 2019, volume 15 issue I,57-70 .

6 .V. E. Tarasov: "Exact Discrete Analogs of Derivatives of Integer Orders: Differences as Infinite Series", Journal of Mathematics, Vol 2015 Article ID 134842.

7. W. Koepf and D.Schmersau (1996) "On the de Branges theorem", Complex Variables, Theory and Application: An International Journal, 31:3, 213-230.

8. D. Merlini, M. Sala, and N. Sala: "A numerical comment on the tiny oscillations and a heuristic conjecture", ArXiv; 1903.623v2 (math NT). (2019).

9. D. Merlini, M. Sala, and N. Sala: "Quasi Fibonacci approximation to the low tiny fluctuations of the Li-Keiper coefficients; a numerical computation", ArXiv :1904.07005 (math GM)(2019).

10. A.Voros: "Sharpening of Li's criterion for the Riemann Hypothesis", ArXiv:math/040421342 (math. N.T.), 2004.

11. Example programs-ARB2.17-gitdocumentation, Li-Coefficients. Retrieved, 10 May 2019: http://arblib.org/examples.html

Appendix 1

The K function i.e. the Koebe function was important in the de Branges proof of the Bieberbach conjecture (de Branges's Theorem). In the variable s it is given by $K(s) = s \cdot (s-1)$; thus in the variable $z=1-1/s$, that is $s=1/(1-z)$ we have:

$$K(z) = z/(1-z)^2 = z + 2 \cdot z^2 + 3 \cdot z^3 + \ldots \sum_{n=0}^{\infty} a_n \cdot z^n, \quad a_0 = 0, a_1 = 1, a_n = n, n \geq 2$$

Notice that for $Re(s) \geq \tfrac{1}{2}$, $s \to z$ maps the half plane $Re(s) \geq \tfrac{1}{2}$ onto the unit circle $|z| \leq 1$. The K function appears naturally in our Eq.(4) i.e. in $f(1/(1-z)) = (1/\gamma) \cdot K(z) \cdot \Psi(z)$ and we have specified $(1/\gamma) \cdot \Psi(z)$ (coefficient of $K(z)$) as a perturbation around the K function which should reduces the coefficients $a_n = n$ to $|b_n| \leq |a_n|$ i.e. $b_n \leq n$.

Appendix 2

Here, we add some computations concerning the possible univalent property of the function of interest, of course, difficult to proove in general and we limit us to look at two particular cases.

1. We show that on the straight line $s = b+i$, $b \in [1, \infty[$, the function f of Eq.(1), is univalent. Below we present the plots of $Re(f(s))$ and $Im(f(s))$ that is real and imaginary part of $f(s)$ on the straight line defined above.

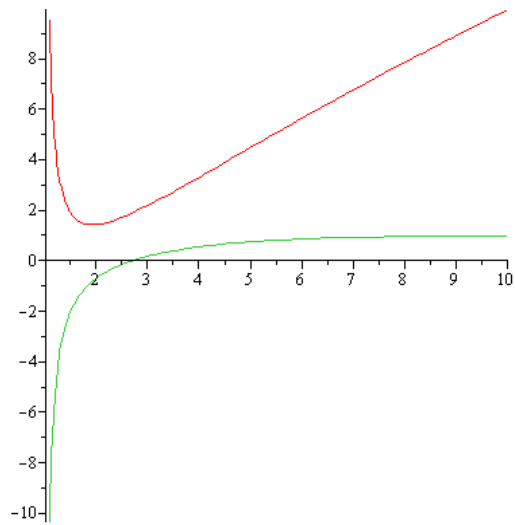

Fig. A.1 In red $Re(f(b+i))$ where $b = Re(s) = Re(b+i \cdot t)$, which is not injective and in green $Im(f(s)) = Im(f(b+i))$ which is injective. Thus $f(s)$ for this case is univalent.

2. A second case is that of a straight line, i.e. b= constant s = 1+i +i.t.

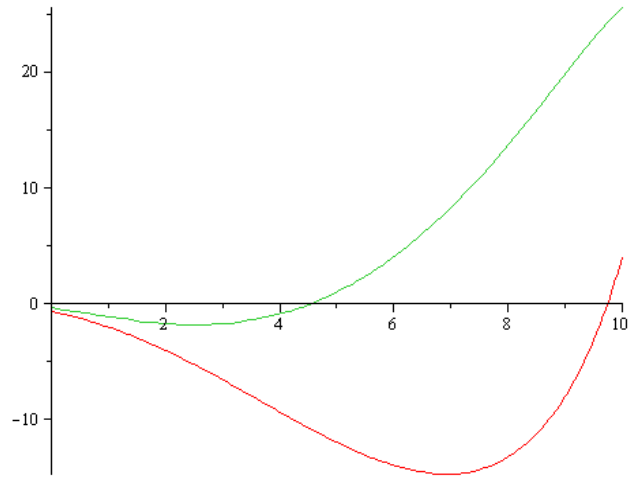

Fig.A.2  In red Re(f(s)), in green Im(f(s)) in the interval t= 0-4.5
Im(f(s)) is not injective ; Re(f(s)) is not injective in 0-9.8
but "injective in 0-4.5. The function is univalent in 0-4.5.
The same in the interval 4.5-∞ where Im(f(s)) is injective.

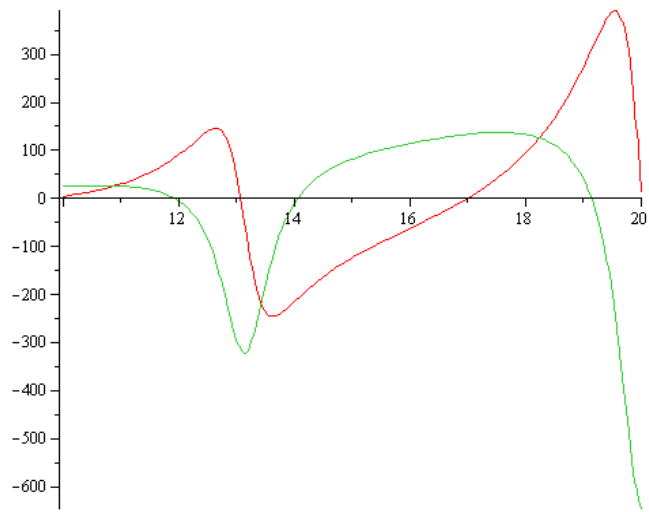

Fig. A.3   In red Re(f(s)), in green Im(f(s)) up to t=20.
f(s) is univalent up to t=20.

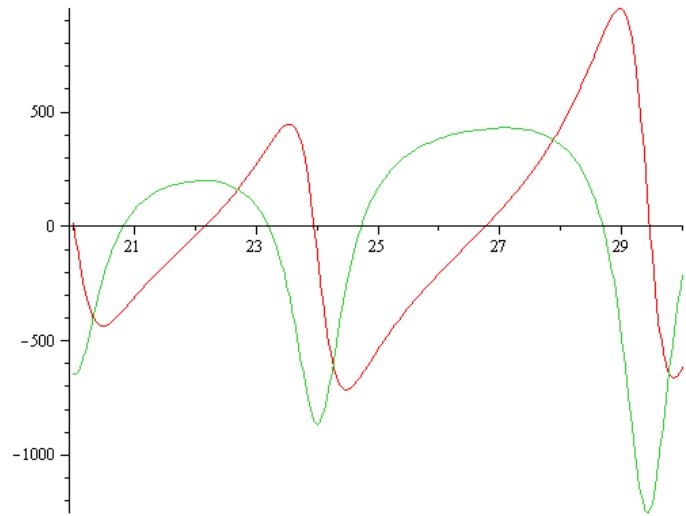

Fig.A.4 In red Re(f(s)), in green Im(f(s)) for t ϵ [20, 30]
f(s) is univalent up to t=30.

The situation is that it is impossible to find 2 values of t such that $Re(f(t_1)) = Re(f(t_2))$ and at the same time $Im(f(t_1)) = Im(f(t_2))$, with $t_1$ different from $t_2$.

## Appendix 3

The system of Equations we have derived without and with the use of the definition of the coefficients have still many solutions and in fact are connected with the discrete derivative of a function as follow. In Ref [6], the finite difference of order n (n integer) for functions of a discrete variable m (in Z) are defined by:

$$\Delta^n f(m) = \sum_{k=0}^{n} (-1)^k \cdot \binom{n}{k} \cdot f(m+(n-k))$$

Setting $\Delta^n f(m) = 0$ and m=0, we have:

$$\sum_{k=0}^{n} (-1)^k \cdot \binom{n}{k} \cdot f(n-k) = 0$$

and thus setting $f(k) = \lambda(k)$ :

$$\lambda(n) = \sum_{k=1}^{n-1} (-1)^{(k+1-n.)} \cdot \binom{n}{k} \cdot \lambda(k)$$

which is our general formula for the complete λ's (valid also for the $\lambda_{tiny}$'s and the $\lambda_{trend}$'s ).